\documentclass[12pt]{amsart}

\usepackage{amsmath}
\usepackage{amsthm}
\usepackage{amssymb}
\usepackage{graphics}
\usepackage{pgf}
\usepackage{tikz}
\usepackage{pifont}
\usetikzlibrary{arrows}
\usepackage[all]{xy}
\usepackage{lineno}

\def\frk{\frak}               

\def\Phi{{\frk n}}
\def\Phi{{\frk N}}
\def\opn#1#2{\def#1{\operatorname{#2}}} 
\opn\chara{char} \opn\length{\ell} \opn\pd{pd} \opn\rk{rk}
\opn\projdim{proj\,dim} \opn\injdim{inj\,dim} \opn\rank{rank}
\opn\depth{depth} \opn\grade{grade} \opn\height{height}
\opn\embdim{emb\,dim} \opn\codim{codim}

\opn\Tr{Tr} \opn\bigrank{big\,rank}
\opn\superheight{superheight}\opn\lcm{lcm}
\opn\trdeg{tr\,deg}
\opn\reg{reg} \opn\lreg{lreg} \opn\ini{in} \opn\lpd{lpd}
\opn\size{size}\opn\bigsize{bigsize}
\opn\cosize{cosize}\opn\bigcosize{bigcosize}
\opn\sdepth{sdepth}\opn\sreg{sreg}
\opn\link{link}\opn\fdepth{fdepth}
\opn\index{index}
\opn\index{index}
\opn\indeg{indeg}
\opn\N{N}
\opn\mult{mult}
\opn\SSC{SSC}
\opn\SC{SC}
\opn\lk{lk}
\opn\HS{HS}
\opn\div{div} \opn\Div{Div} \opn\cl{cl} \opn\Cl{Cl}
\opn\Spec{Spec} \opn\Supp{Supp} \opn\supp{supp} \opn\Sing{Sing}
\opn\Ass{Ass} \opn\Min{Min}\opn\Mon{Mon} \opn\dstab{dstab} \opn\astab{astab}
\opn\Syz{Syz}
\opn\reg{reg}
\opn\Ann{Ann} \opn\Rad{Rad} \opn\Soc{Soc}
\opn\Im{Im} \opn\Ker{Ker} \opn\Coker{Coker} \opn\Am{Am}
\opn\Hom{Hom} \opn\Tor{Tor} \opn\Ext{Ext} \opn\End{End}\opn\Der{Der}
\opn\Aut{Aut} \opn\id{id}

\opn\nat{nat}
\opn\pff{pf}
\opn\Pf{Pf} \opn\GL{GL} \opn\SL{SL} \opn\mod{mod} \opn\ord{ord}
\opn\Gin{Gin} \opn\Hilb{Hilb}\opn\sort{sort}
\opn\initial{init}
\opn\ende{end}
\opn\height{height}
\opn\type{type}
\opn\aff{aff} \opn\con{conv} \opn\relint{relint} \opn\st{st}
\opn\lk{lk} \opn\cn{cn} \opn\core{core} \opn\vol{vol}
\opn\link{link} \opn\Link{Link}\opn\lex{lex}
\opn\gr{gr}

\def\pot#1#2{#1[\kern-0.28ex[#2]\kern-0.28ex]}

\opn\dirlim{\underrightarrow{\lim}}
\opn\inivlim{\underleftarrow{\lim}}
\let\To=\longrightarrow
\def\Implies{\ifmmode\Longrightarrow \else
        \unskip${}\Longrightarrow{}$\ignorespaces\fi}
\def\implies{\ifmmode\Rightarrow \else
        \unskip${}\Rightarrow{}$\ignorespaces\fi}
\def\iff{\ifmmode\Longleftrightarrow \else
        \unskip${}\Longleftrightarrow{}$\ignorespaces\fi}

\let\:=\colon
\newtheorem{Theorem}{Theorem}[section]
 \newtheorem{Lemma}[Theorem]{Lemma}
 \newtheorem{Corollary}[Theorem]{Corollary}
 
 \newtheorem{Remark}[Theorem]{Remark}
 
 \newtheorem{Example}[Theorem]{Example}
 
 \newtheorem{Definition}[Theorem]{Definition}

\newtheorem{Notation}[Theorem]{Notation}

\let\epsilon\varepsilon
\let\kappa=\varkappa
\def\qed{\ifhmode\textqed\fi
      \ifmmode\ifinner\quad\qedsymbol\else\dispqed\fi\fi}
\def\textqed{\unskip\nobreak\penalty50
       \hskip2em\hbox{}\nobreak\hfil\qedsymbol
       \parfillskip=0pt \finalhyphendemerits=0}
\def\dispqed{\rlap{\qquad\qedsymbol}}

\opn\dis{dis}
\def\pnt{{\raise0.5mm\hbox{\large\bf.}}}

\opn\Lex{Lex}

\opn\Spec{Spec} \opn\Supp{Supp} \opn\supp{supp}
 \opn\Ass{Ass}
 \opn\p{Ass}
 \opn\min{min}
 \opn\max{max}
 \opn\MIN{Min}
 \opn\p{\mathfrak{p}}
\opn\Deg{Deg}

\begin{document}

 \title{The Multiplicity of Powers of a Class of Non-Square-free Monomial Ideals}

\author {Liuqing Yang, Zexin Wang*}

\footnote{* Corresponding author}
\address{School  of Mathematical Science, Soochow University, 215006 Suzhou, P.R.China}
\email{zexinwang6@outlook.com}

\address{School of Mathematical Science,  Soochow University, 215006 Suzhou, P.R.China}
\email{20214007001@stu.suda.edu.cn}

 \begin{abstract} Let $R = \mathbb{K}[x_1, \ldots, x_n]$ be a polynomial ring over a field $\mathbb{K}$, and let $I \subseteq R$ be a monomial ideal of height $h$. We provide a formula for the multiplicity of the powers of $I$ when all the primary ideals of height $h$ in the irredundant reduced primary decomposition of $I$ are irreducible. This is a generalization of \cite[Theorem 1.1]{TV}. Furthermore, we present a formula for the multiplicity of powers of special powers of monomial ideals that satisfy the aforementioned conditions. Here, for an integer $m>0$, the $m$-th special power of a monomial ideal refers to the ideal generated by the $m$-th powers of all its minimal generators. Finally, we explicitly provide a formula for the multiplicity of powers of special powers of edge ideals of weighted oriented graphs.
 \end{abstract}



\maketitle

\section{Introduction}

 Let $R=\mathbb{K}[x_1,\ldots,x_n]$ be a standardly graded   polynomial ring over a field $\mathbb{K}$,  and let $M$ be a finitely generated graded $R$-module. We denote by $M_k$  the degree $k$ component of $M$. The {\it Hilbert function} $H_M(k)$ of $M$ is a function from $\mathbb{Z}$ to $\mathbb{N}$ given by $H_M(k):=\dim_{\mathbb{K}}M_k$   for all  $k\in \mathbb{ Z}$. The {\it Hilbert series} $\HS(M,t)$ of $M$ is defined to be the formal power series: $$\HS(M,t):=\sum_{k\in \mathbb{Z}}H_M(k)t^k.$$
Assuming $\dim M=d+1$, David Hilbert showed that  that $\HS(M,t)$ is a rational function of the following form $$\HS(M,t)=\frac{Q_M(t)}{(1-t)^{d+1}}.$$
Here, $Q_M(t)\in \mathbb{Q}[t,t^{-1}]$ is  a Laurent polynomial such that $Q_M(1)\neq 0$. As  a consequence,  there exists a polynomial $p_M(t)\in \mathbb{Q}[t]$ of degree $d$ such that $H_M(k)=p_M(k)$ for all $k\gg 0$. The polynomial $p_M(t)$ is referred to as the {\it Hilbert polynomial} of $M$.

\begin{Definition} \em Let $M$ be a finitely generated graded $R$-module of dimension $d+1$. The  Hilbert polynomial $p_M(t)$ of $M$ can be written as $$p_M(t)=\sum_{i=0}^d(-1)^ie_i(M)\binom{t+d-i}{d-i}.$$
The integer coefficients $e_i(M)$ for $i=0,\ldots,d$ are called the {\it Hilbert coefficients} of $M$.
\end{Definition}

 According to \cite[Proposition 4.1.9]{BH}, we have $e_i(M)=\frac{Q_M^{(i)}(1)}{i!}$ for $i=0,\ldots,d$. The first Hilbert coefficients $e_0(M)$ is also called the {\it multiplicity} of $M$ and denoted by $\mult(M)$. The multiplicity of a graded ideal is a significant  invariant in algebraic geometry and commutative algebra. For the  study of the multiplicity of graded ideals, we refer to \cite{A1,A2}, \cite{HS1,HS2}, \cite{TY1,TY2,TY3} and the references therein.

 Let $I$ be a graded ideal   with $\dim R/I=d$. Herzog-Puthenpurakal-Verma  showed in \cite[Theorem 1.1]{HPV08} that {\it $e_i(R/I^s)$ is of polynomial type in $s$ of degree $\leq n-d+i$ for $i=0,1,\ldots,d$}. Recall that a function $f:\mathbb{N}\rightarrow \mathbb{Q}$ is of {\it polynomial type} of degree $d$ if there exists a polynomial $p(t)\in \mathbb{Q}[t]$ of degree $d$ such that $f(k)=p(k)$ for all $k\gg 0$.
 In particular, $\mult(R/I^s)$ is of polynomial type in $s$. Naturally, the question arises: Is it possible to explicitly compute the multiplicity of powers of graded ideals in certain instances?

To the best of our knowledge, the first instance of a non-trivial graded ideal, for which the multiplicity of its powers is explicitly provided, is the path ideal of a line graph, as computed previously in \cite{SWL}, where it was proved that if $I$ is the path ideal of a line graph, then the formula  $\mult(R/I^s) = \mult(R/I) \binom{\height(I)+s-1}{s-1}$ holds. Recently, Thuy and Vu extended this formula to encompass arbitrary square-free monomial ideals in their work \cite{TV}. Building upon their findings, we further extend \cite[Theorem 1.1 and Lemma 2.1]{TV} to compute the multiplicity of powers for a specific class of non-square-free monomial ideals, notably including the edge ideals of weighted oriented graphs and edge-weighted graphs.
\begin{Definition}\em
    Let $I$ be a graded ideal of $R$ with a primary decomposition given by
    \[
    I = Q_1 \cap Q_2 \cdots \cap Q_t.
    \]
    We refer to this decomposition as \emph{reduced} if the radicals $\sqrt{Q_i}$ are pairwise distinct for all $i = 1, \ldots, t$. Furthermore, we refer to it as \emph{irredundant} if for any $1 \leq i \leq t$, the ideal $Q_i$ is not a superset of the intersection of the other primary ideals, i.e., $Q_i \nsupseteq \bigcap_{j \neq i} Q_j$.
\end{Definition}

\begin{Theorem}{\em (Theorem~\ref{Main result})}
Let $I$ be a monomial ideal of $R$ with height $h$. Suppose $I$ admits an irredundant reduced primary decomposition
\[
I = Q_1 \cap \cdots \cap Q_r \cap Q_{r+1} \cap \cdots \cap Q_t,
\]
where $\height(Q_i) = h$ for $i = 1,\ldots,r$ and $\height(Q_i) > h$ for $i = r+1,\ldots,t$. Then following statements hold.

\textup{(1)} For every integer $s \geq 1$,
\[
\mult(R/I^s) = \sum_{i=1}^r \mult(R/Q_i^s).
\]

\textup{(2)} If each $Q_i$ ($1 \leq i \leq r$) is an irreducible monomial ideal generated by pure powers of variables with exponents $a_{i_1},\ldots,a_{i_h}$, then for any $s \geq 1$,
\[
\mult(R/I^s) = \mult(R/I) \binom{h+s-1}{s-1} = \sum_{i=1}^r \left( \prod_{j=1}^h a_{i_j} \right) \binom{h+s-1}{s-1}.
\]
\end{Theorem}

We remark that \cite[Lemma 2.1]{TV} is a special case of formula (\dag) when $s=1$, and \cite[Theorem 1.1]{TV} is a special case of formula (\ddag) when $I$ is a square-free monomial ideal.

Let $I$ be a monomial ideal. The ideal generated by the $m$-th powers of all its minimal generators is called the $m$-th \emph{special power} of $I$, denoted as $I^{\{m\}}$. We have proved the following theorem:

\begin{Theorem} {\em (Theorem~\ref{main2})}
If $I$ satisfies the hypotheses of Theorem~\ref{Main result}(2), then $I^{\{m\}}$ also satisfies them for all integers $m \geq 1$.
Furthermore, let $\height(I)=h$, then for all $m,s \geq 1$,
\[
\mult(R/(I^{\{m\}})^s) = m^h \begin{pmatrix}h+s-1\\ s-1\end{pmatrix} \mult(R/I).
\]
\end{Theorem}

We provide some notations and definitions that will be used throughout this paper.
\begin{Notation}\label{graph} \em
Let $G = (V(G), E(G))$ be a simple graph (without loops or multiple edges) with vertices $V(G) = \{x_1, \ldots, x_n\}$ and edge set $E(G)$. By identifying the variables of the polynomial ring $R = \mathbb{K}[x_1, \ldots, x_n]$ with the vertices of $V(G)$, we can associate to $G$ a square-free monomial ideal $I(G) = (\{x_ix_j \mid \{x_i, x_j\} \in E(G)\})$, called the edge ideal of $G$. 
\end{Notation}
\begin{Definition}\label{cover} \em  For a vertex $x_i\in V(G)$, the {\it neighbor set} of $x_i$ is defined to be the set $N_G(x_i) = \{x_j |\ \{x_i,x_j\}\in E(G)\} $. A {\it vertex cover} of $G$ is a subset $C \subseteq V(G)$ such that for each edge $\{x_i, x_j\}$ in $G$, either $x_i \in C$ or $x_j \in C$. A vertex cover is {\it minimal} if it does not properly contain another vertex cover of $G$. The minimum number of vertices in a minimal vertex cover of $G$ is called the {\it vertex covering number} of $G$, denoted as $\alpha(G)$. Let $r(G)$ denote the number of minimal vertex covers of $G$ that contain exactly $\alpha(G)$ vertices.
\end{Definition}

\begin{Definition}\label{oriented graph} \em  A {\it weighted oriented graph} $D$, whose underlying graph is $G$, is a triplet $(V(D), E(D), w)$ where $V(D) = V(G)$, $E(D) \subseteq V(D)\times V(D)$ such that
 $\{\{ x_i, x_j \} | (x_i, x_j) \in E(D)\} = E(G)$, and $w$ is a function $w: V(D) \rightarrow \mathbb{N}$. The vertex set of $D$ and the edge set of $D$ are $V(D)$ and $E(D)$, respectively.
Sometimes, for brevity, we denote these sets by $V$ and $E$ respectively.
 The {\it weight} of $x_i \in V$ is $w(x_i)$.
\end{Definition}

\begin{Definition}\label{L} \em The edge ideal of a weighted oriented graph $D$ is a monomial ideal given by
\[
I(D) = (x_i x_j^{w(x_j)} \mid (x_i, x_j) \in E(D)) \subseteq R.
\]
\end{Definition}
Edge ideals of weighted oriented graph arose in the theory of Reed-Muller codes as initial
ideals of vanishing ideals of projective spaces over finite fields (see \cite{MPV}, \cite{PS}). In recent years, its algebraic properties have been studied by many researchers. Relevant research can be referred to in  \cite{CK}, \cite{HLMRV}, \cite{MP}, \cite{WZXZ} and \cite{ZXWT}, etc. We provide a formula for the multiplicity of powers of special powers of the edge ideal of any weighted oriented graph using combinatorial properties.
\begin{Definition} \em Let $D = (V, E, w)$ be a vertex-weighted  oriented graph and $G$ its underlying graph. For a vertex cover $C$ of $G$, define
\begin{align*}
&L_1(C) = \{x_i \in C\  | \  \exists \ x_j \text{ such that } (x_i,x_j) \in E \text{ and } x_j \notin C\}, \\
&L_3(C) = \{x_i \in C\  | \ N_G(x_i) \subseteq C\}, \\
&L_2(C) = C \setminus (L_1(C) \cup L_3(C)).
\end{align*}

\end{Definition}

\begin{Theorem}{\em (Theorem~\ref{main3})} \em Let $D = (V, E, w)$ be a vertex-weighted  oriented graph and $G$ its underlying graph. Let  $C_1, \ldots, C_{r(G)}$ be all minimal vertex covers of graph $G$ that contain exactly  $\alpha(G)$ vertices.  Then, for all $m, s \geq 1$,
   $$\mult(R/(I(D)^{\{m\}})^s) = m^{\alpha(G)} \sum_{i=1}^{r(G)} \left( \prod_{x_j\in L_2(C_i)} w(x_j)\right) \binom{\alpha(G)+s-1}{s-1} .$$
\end{Theorem}
When $m=1$, this formula reduces to the multiplicity formula for powers of the edge ideal of a weighted oriented graph.

The paper is structured as follows: Section 2 provides an explicit formula for the multiplicity of the powers of $I$ when all the primary ideals of height $h$ in the irredundant reduced primary decomposition of $I$ are irreducible. Then, we introduce the concept of special powers, and further derive a formula for the multiplicity of the powers of these special powers for such ideals. Section 3 provides a formula for the multiplicity of powers of special powers of the edge ideal of any weighted oriented graph using combinatorial properties.

\section{Multiplicity}

In this section, we will present a formula for computing the multiplicity of the powers of a class of non-square-free monomial ideals. Then, we introduce the concept of special powers, and further derive a formula for the multiplicity of the powers of these special powers for such ideals. Throughout this paper, the polynomial ring $\mathbb{K}[x_1, \ldots, x_n]$ will be uniformly denoted by $R$.
\begin{Definition}\em
    Let $I$ be a graded ideal of $R$ with a primary decomposition given by
    \[
    I = Q_1 \cap Q_2 \cdots \cap Q_t.
    \]
    We refer to this decomposition as \emph{reduced} if the radicals $\sqrt{Q_i}$ are pairwise distinct for all $i = 1, \ldots, t$. Furthermore, we refer to it as \emph{irredundant} if for any $1 \leq i \leq t$, the ideal $Q_i$ is not a superset of the intersection of the other primary ideals, i.e., $Q_i \nsupseteq \bigcap_{j \neq i} Q_j$.
\end{Definition}

\begin{Theorem}\label{Main result}
Let $I$ be a monomial ideal of $R$ with height $h$. Suppose $I$ admits an irredundant reduced primary decomposition
\[
I = Q_1 \cap \cdots \cap Q_r \cap Q_{r+1} \cap \cdots \cap Q_t,
\]
where $\height(Q_i) = h$ for $i = 1,\ldots,r$ and $\height(Q_i) > h$ for $i = r+1,\ldots,t$. Then following statements hold.

\textup{(1)} For every integer $s \geq 1$, we have
\[
\mult(R/I^s) = \sum_{i=1}^r \mult(R/Q_i^s).
\]

\textup{(2)} If each $Q_i$ ($1 \leq i \leq r$) is an irreducible monomial ideal generated by pure powers of variables with exponents $a_{i_1},\ldots,a_{i_h}$, then for any $s \geq 1$, we have
\[\tag{\dag}
\mult(R/I^s) = \mult(R/I) \binom{h+s-1}{s-1} = \sum_{i=1}^r \left( \prod_{j=1}^h a_{i_j} \right) \binom{h+s-1}{s-1}.
\]
\end{Theorem}

We first prove that  formula $(\dag)$ holds for irreducible monomial ideals.
According to \cite[Corollary 1.3.2]{HH}, a monomial ideal is irreducible if and only if it is generated by pure powers of variables. As usual, we use $G(I)$ to denote the minimal generating set for a monomial ideal $I$.

\begin{Notation}\label{Notation}
 Without loss of generality, let $I_m = (x_1^{a_1}, x_2^{a_2}, \ldots, x_m^{a_m})$ be an irreducible monomial ideal of $R$, where $1 \leq m \leq n$ and $a_i$ are positive integers for all $i$. When $m > 1$, denote $I_{m-1} = (x_1^{a_1}, x_2^{a_2}, \ldots, x_{m-1}^{a_{m-1}})$.
 \end{Notation}

\begin{Lemma} \label{colon}
 Following Notation \ref{Notation}, for every integer $s\geq 2$,
 $$I_m^s:x_m^{a_m}=I_m^{s-1}.$$
\end{Lemma}
   \begin{proof}
Because $x_m^{a_m} \in G(I_m)$, the  relation ``$\supseteq$'' is obvious.
 We only need to prove the reverse inclusion ``$\subseteq$'' holds. Let $u \in I_m^s:x_m^{a_m}$ be a monomial, then $ux_m^{a_m} \in I_m^s$. This means there exists $v \in G(I_m^s)$ such that $v | ux_m^{a_m}$. Note that we may write $v$ as $x_1^{k_1a_1}x_2^{k_2a_2}\cdots x_m^{k_ma_m}$, where $k_1+\cdots+k_m=s$ with $k_i\geq 0$ for $i=1,\ldots,m$.
If $k_m=0$, then $v | u$, and thus $u \in I_m^s \subseteq I_m^{s-1}$. If $k_m\geq 1$, then, since $\frac{v}{x_m^{a_m}}=x_1^{k_1a_1}x_2^{k_2a_2}\cdots x_m^{(k_m-1)a_m}$, we conclude that $\frac{v}{x_m^{a_m}}$ belongs to $I_m^{s-1}$. It follows that $u$ belongs to $I_m^{s-1}$ since $\frac{v}{x_m^{a_m}} \mid u$, as required.
   \end{proof}

   \begin{Lemma} \label{irreducible}
     Following Notation \ref{Notation}, for every integer $s\geq 1$,  $$\mult (R/I_m^s)=\mult (R/I_m)\begin{pmatrix}m+s-1\\s-1\end{pmatrix}=a_1\ldots a_m\begin{pmatrix}m+s-1\\s-1\end{pmatrix}.$$
\end{Lemma}
\begin{proof} Since $I_m$ is a monomial ideal generated by a regular sequence of monomials, according to \cite[Exercise~16.9]{P}, we have $\mult(R/I_m) = a_1\ldots a_m$. Therefore, it suffices to prove that the first term  is equal to the third term in the equality.

We proceed by induction on both $s$ and $m$.
      The case when $s=1$ or $m=1$  follows from \cite[Exercise~16.9]{P}.      Suppose now that $s>1$ and  $m>1$.
 Consider the following two short exact sequences of graded $R$-modules:

\begin{equation}
     0 \To\frac{R}{I_m^s:x_m^{a_m}}[-a_m] \To\frac{R}{I_m^s} \To\frac{R}{(I_m^s,x_m^{a_m})} \To 0,
     \end{equation}

 \begin{equation}
   0 \To\frac{R}{I_{m-1}^{s}:x_m^{a_m}}[-a_m] \To\frac{R}{I_{m-1}^{s}} \To\frac{R}{(I_{m-1}^{s},x_m^{a_m})} \To 0.
     \end{equation}

  Firstly, by Lemma~\ref{colon}, we have $\sqrt{I_m^s:x_m^{a_m}} = \sqrt{I_m^{s-1}} = (x_1, \ldots, x_m)$. Therefore,
\[
\dim  \frac{R}{I_m^s:x_m^{a_m}} = \dim  \frac{R}{I_m^{s-1}} = n - m.
\]

Note that $(I_m^s,x_m^{a_m})=(I_{m-1}^s,x_m^{a_m})$ and $\sqrt{(I_{m-1}^s,x_m^{a_m})}=(x_1,\ldots,x_m)$,  we obtain  $$\dim \frac{R}{(I_m^s,x_m^{a_m})}=\dim \frac{R}{(I_{m-1}^s,x_m^{a_m})}=n-m.$$

Applying \cite[Lemma 3.9]{SWL} to the exact sequence (1), we obtain
\[
\mult(R/I_m^s) = \mult(R/I_m^{s-1}) + \mult(R/(I_{m-1}^s, x_m^{a_m})).
\]
 Next, consider the second short exact sequence, where $I_{m-1}^{s}:x_m^{a_m}=I_{m-1}^{s}$ holds trivially, leading to the following equality:
\[
 \HS(R/(I_{m-1}^{s},x_m^{a_m}),t) =(1-t^{a_m})\HS(R/I_{m-1}^s,t)=\frac{(1+t+\cdots+t^{a_m-1})Q(t)}{(1-t)^{d-1}}.
\]
Here, we assume that $ \dim \frac{R}{I_{m-1}^s}=d$ and $\HS(R/I_{m-1}^s,t)=\frac{Q(t)}{(1-t)^d}$.
Since $\sqrt{I_{m-1}^s} = (x_1, \ldots, x_{m-1})$ and $\sqrt{(I_{m-1}^{s}, x_m^{a_m})} = (x_1, \ldots, x_m)$, we obtain
\[
\dim \frac{R}{(I_{m-1}^s, x_m^{a_m})} = \dim \frac{R}{I_{m-1}^s} - 1=d-1.
\]
It follows  \cite[Proposition 4.1.9]{BH} that
\[
\mult(R/(I_{m-1}^s, x_m^{a_m})) = a_m \mult(R/I_{m-1}^s).
\]
Combining the above equality with  the induction hypothesis,  we can deduce that:
    \begin{align*}
      & \mult (R/I_m^s)= \mult (R/I_m^{s-1})+a_m\mult (R/I_{m-1}^s)
      \\
            =&a_1\ldots a_m\begin{pmatrix}m+s-2\\ s-2\end{pmatrix}
      +a_m a_1\ldots a_{m-1}\begin{pmatrix}m+s-2\\s-1\end{pmatrix}\\
            =&a_1\ldots a_m\begin{pmatrix}m+s-1\\s-1\end{pmatrix}.
    \end{align*}
    \end{proof}

\begin{proof}[Proof of Theorem \ref{Main result}]
Suppose that $I^s$ admits an
irredundant reduced primary decomposition
\[
I^s = Q_1' \cap \cdots \cap Q_r' \cap Q_{r+1}' \cap \cdots \cap Q_k',
\]
such that $\sqrt{Q_i'}$ are pairwise distinct for $i=1,\ldots,k$. Assume further that $\sqrt{Q_i'} = \sqrt{Q_i}$ for $i = 1, \ldots, r$.
Since $\sqrt{Q_1'}, \ldots, \sqrt{Q_r'}$ are all minimal prime ideals of $I$,  it is easy to see that
\[
Q_i' = Q_i^s, \text{ for any } 1 \leq i \leq r.
\]
One may also look at the proof of \cite[Lemma 2]{GMSVV} for this observation.
So, applying \cite[Lemma 2.1]{TV} directly yields the first assertion.
The second assertion follows from Theorem \ref{Main result} (1) and Lemma~\ref{irreducible}.
\end{proof}
Using the concept of special powers of monomial ideals as defined below, we can construct numerous monomial ideals satisfying the hypotheses of Theorem~\ref{Main result}(2).
\begin{Definition} \em Let $I$ be a monomial ideal with $G(I)=\{u_1,\ldots,u_t\}$. For any $m \geq 1$, define the $m$-th {\it special power} of $I$ as
\[
I^{\{m\}} = (u_1^m, \ldots, u_t^m).
\]
\end{Definition}

We collect some easy facts regarding the special power.

\begin{Lemma}\label{special power basic properties}
Let $I_1, I_2, \ldots, I_t$ be monomial ideals, then for $m \geq 1$, we have:
\begin{enumerate}
 \item $(I_1\cap I_2\cap \cdots \cap I_t)^{\{m\}}=I_1^{\{m\}}\cap I_2^{\{m\}}\cap \cdots \cap I_t^{\{m\}}$;
  \item $(I_1I_2\cdots I_t)^{\{m\}}=I_1^{\{m\}}I_2^{\{m\}}\cdots I_t^{\{m\}}$;
  \item $(I_1+I_2+\cdots +I_t)^{\{m\}}=I_1^{\{m\}}+I_2^{\{m\}}+\cdots +I_t^{\{m\}}$;
  \item $I$ is an irreducible monomial ideal if and only if so is $I^{\{m\}}$;
   \item $I$ is a primary  monomial ideal if and only if so is $I^{\{m\}}$.
  \end{enumerate}
\end{Lemma}
\begin{proof} (1)  By induction, we only consider the case $t=2$. For any $u \in G(I_1 \cap I_2)$, it holds that $u \in I_1,I_2$, and thus $u^m \in I_1^{\{m\}}, I_2^{\{m\}}$ for all $m$. Since all $u^m$ generate $(I_1 \cap I_2)^{\{m\}}$, we have $(I_1 \cap I_2)^{\{m\}} \subseteq I_j^{\{m\}}$ for $j=1,2$, which further implies $(I_1 \cap I_2)^{\{m\}} \subseteq  I_1^{\{m\}}\cap I_2^{\{m\}}$.
Next, we prove the reverse inclusion. For any $u \in G(I_1)$ and any $v \in G(I_2)$, we have $\text{lcm}(u^m, v^m) = (\text{lcm}(u, v))^m \in (I_1 \cap I_2)^{\{m\}}$ for all $m$. And since all $\text{lcm}(u^m, v^m)$ generate $I_1^{\{m\}} \cap I_2^{\{m\}}$, the reverse inclusion also holds.

The proofs of (2) and (3) is similar as (1) and we omit it.

(4) The conclusion follows immediately from the generating structure of irreducible monomial ideals by pure powers of variables.

(5) A monomial is a primary ideal if and only if it is the intersection of irreducible monomial ideals with the same support. In view of this fact, the assertion follows from (1) together with (4).
\end{proof}

\begin{Lemma}\label{special power property}
If a monomial ideal $I$ admits an irredundant reduced primary decomposition $I = \bigcap_{i=1}^t Q_i$, then
\[
I^{\{m\}} = \bigcap_{i=1}^t Q_i^{\{m\}}
\]
is an irredundant reduced primary decomposition of $I^{\{m\}}$.

\end{Lemma}

\begin{proof}
From Lemma~\ref{special power basic properties}(1) and (5), we can deduce that
\[\tag{\S}
I^{\{m\}} = \bigcap_{i=1}^t Q_i^{\{m\}}
\]
is a primary decomposition of $I^{\{m\}}$. Furthermore, since $\sqrt{Q_i} = \sqrt{Q_i^{\{m\}}}$, and given that $I = \bigcap_{i=1}^t Q_i$ is a reduced primary decomposition, it can be concluded that $(\S)$ is also a reduced primary decomposition. We only need to prove that this decomposition is irredundant.

If there exists an $i$ such that $Q_i^{\{m\}} \supseteq \bigcap_{j \neq i} Q_j^{\{m\}}$, then by Lemma~\ref{special power basic properties}(1), we have $Q_i^{\{m\}} \supseteq \left( \bigcap_{j \neq i} Q_j \right)^{\{m\}}$. For any $u \in G\left( \bigcap_{j \neq i} Q_j \right)$, there exists a $v \in G(Q_i)$ such that $v^m | u^m$, thereby $v | u$. Therefore, $Q_i \supseteq \bigcap_{j \neq i} Q_j$, which contradicts the fact that $I = \bigcap_{i=1}^t Q_i$ is irredundant. This completes the proof.
\end{proof}

From Theorem \ref{Main result}, we can derive the following results.
\begin{Theorem}\label{main2}
If $I$ satisfies the hypotheses of Theorem~\ref{Main result}(2), then $I^{\{m\}}$ also satisfies them for all integers $m \geq 1$.
Furthermore, let $\height(I)=h$, then for all $m,s \geq 1$,
\[
\mult(R/(I^{\{m\}})^s) = m^h \begin{pmatrix}h+s-1\\ s-1\end{pmatrix} \mult(R/I).
\]

\end{Theorem}
\begin{proof}

The first assertion directly follows from Lemmas \ref{special power property} and \ref{special power basic properties}(4).
For the second assertion, by Theorem~\ref{Main result}(2), it suffices to prove that
    \[
    \mult(R/I^{\{m\}}) = m^h \mult(R/I).
    \]
If the irreducible monomial ideal $Q$ is generated by pure powers of variables with degrees $a_1, \ldots, a_h$, then $Q^{\{m\}}$ is generated by pure powers of variables with degrees $ma_1, \ldots, ma_h$.  Additionally, according  to \cite[Exercise~16.9]{P}, we have $$\operatorname{mult}(R/Q^{\{m\}}) = m^h a_1 \ldots a_h=m^h \cdot \operatorname{mult}(R/Q).$$ Thus, combining Theorem \ref{Main result}(2) and Lemma \ref{special power property}, the conclusion is obvious.
\end{proof}
Let $I$ be a square-free monomial ideal in $R$. According to \cite[Corollary~6.2.3]{HH}, the multiplicity $\text{mult}(R/I)$ is the count of associated prime ideals of $I$ having the minimal height. Furthermore, since a square-free monomial ideal can be decomposed as the intersection of monomial prime ideals, by applying Theorem \ref{main2}, we can derive the following corollary.

\begin{Corollary} \label{special power square-free}
\em Let $I$ be a square-free monomial ideal of height $h$ in $R$, and let $r$ be the number of height-$h$ associated primes of $R/I$. Then, for all $m, s \geq 1$,
   $$\mult(R/(I^{\{m\}})^s) = r m^h\begin{pmatrix}h+s-1\\ s-1\end{pmatrix}.$$
\end{Corollary}

\begin{Remark}\em When $m=1$, the above corollary reduces to \cite[Theorem 1.1]{TV}.
\end{Remark}

In the conclusions below, we use the notation from Notation \ref{graph} and Definition \ref{cover}.

\begin{Corollary}\label{special power graph} \em Let $I(G)$ be the edge ideal of the graph $G$. Then, for all $m, s \geq 1$,
   $$\mult(R/(I(G)^{\{m\}})^s) = r(G) m^{\alpha(G)}\begin{pmatrix}\alpha(G)+s-1\\ s-1\end{pmatrix}.$$
\end{Corollary}
\begin{proof}
We know that $\height(I(G))=\alpha(G)$ and $\mult(R/I(G))= r(G)$, so the conclusion is a direct corollary of Corollary \ref{special power square-free}.
\end{proof}
\begin{Remark}\em The $m$-th special power $I(G)^{\{m\}}$ of $I(G)$ is precisely the edge ideal of trivially edge-weighted graphs $G_w$ where each edge has a weight of $m$, as defined in \cite{PS}.

\end{Remark}

\section{weighted oriented graph}

To compute multiplicities of powers for special powers in edge ideals of weighted oriented graphs, we employ Theorem~\ref{Main result} in this section. In this section, we continue to use Notation \ref{graph}, Definition \ref{cover} and Definition \ref{oriented graph}. Recall a  weighted oriented graph $D$, whose underlying graph is $G$, is a triplet $(V(D), E(D), w)$ where $V(D) = V(G)$, $E(D) \subseteq V(D)\times V(D)$ such that
 $\{\{ x_i, x_j \} | (x_i, x_j) \in E(D)\} = E(G)$, and $w$ is a function $w: V(D) \rightarrow \mathbb{N}$. The vertex set of $D$ and the edge set of $D$ are $V(D)$ and $E(D)$, respectively. The edge ideal of a weighted oriented graph $D$ is a monomial ideal given by $I(D) = (x_i x_j^{w(x_j)} \mid (x_i, x_j) \in E(D)) \subseteq R$.

For reference, we repeat Definition \ref{L} as follows:
\begin{Definition}\em Let $D = (V, E, w)$ be a weighted oriented graph and $G$ its underlying graph. For a vertex cover $C$ of $G$, define
\begin{align*}
&L_1(C) = \{x_i \in C\  | \  \exists \ x_j \text{ such that } \{x_i,x_j\} \in E \text{ and } x_j \notin C\}, \\
&L_3(C) = \{x_i \in C\  | \ N_G(x_i) \subseteq C\}, \\
&L_2(C) = C \setminus (L_1(C) \cup L_3(C)).
\end{align*}
\end{Definition}

A vertex cover $C$ of $G$ is called a {\it strong vertex cover} of $D$ if $C$ is a minimal vertex cover of $G$ or for each $x_i \in L_3(C)$ there is an edge $(x_j,x_i) \in E$ such that $x_j \in L_2(C) \cup L_3(C)$ with $w(x_j) \geq 2$.

\begin{Theorem}\label{main3} \em Let $D = (V, E, w)$ be a weighted oriented graph and $G$ its underlying graph. Let  $C_1, \ldots, C_{r(G)}$ be all minimal vertex covers of graph $G$ that contain exactly  $\alpha(G)$ vertices.  Then, for all $m, s \geq 1$,
   $$\mult(R/(I(D)^{\{m\}})^s) = m^{\alpha(G)} \sum_{i=1}^{r(G)} \left( \prod_{x_j\in L_2(C_i)} w(x_j)\right) \binom{\alpha(G)+s-1}{s-1} .$$
\end{Theorem}

\begin{proof}
According to \cite[Remark 26]{PRT}, if $\mathcal{C}_{s}$ denotes the set of strong vertex covers of $D$, then the irredundant reduced primary decomposition of the ideal $I(D)$ is expressed as
\[
I(D) = \bigcap_{C \in \mathcal{C}_{s}} I_{C},
\]
where
\[
I_{C} = \left( L_1(C) \cup \{x_j^{w(x_j)} \mid x_j \in L_2(C) \cup L_3(C)\} \right).
\]
Therefore, $I(D)$ satisfies the hypotheses of Theorem~\ref{Main result}(2).  Using Theorem \ref{main2}, we only need to prove that
\[
\mult(R/I(D)) = \sum_{i=1}^{r(G)} \left( \prod_{x_j \in L_2(C_i)} w(x_j) \right).
\]
According to \cite[Proposition 6]{PRT}, if $C$ is a minimal vertex cover of $D$, then $L_{3}(C) = \emptyset$.
Therefore, for any $1 \leq i \leq r(G)$, we have
$\mult(R/I_{C_i}) = \prod_{x_j \in L_2(C_i)} w(x_j)$. By Theorem~\ref{Main result}(1), the above equality follows as required.
\end{proof}

\begin{Remark} \em From this paper, we know that $\mult(R/I^s) = \mult(R/I) \binom{\height(I)+s-1}{s-1}$ holds for many monomial ideals $I$. Naturally, one would ask whether this rule holds for all monomial ideals. The answer is no, as we will see in the following example.
\end{Remark}
\begin{Example} \em
Let $I = (x_1^2, x_2^2, x_3^4) \cap (x_1^3, x_2^3, x_3^2) = (x_2^3, x_1^3, x_2^2x_3^2, x_1^2x_3^2, x_3^4)$ be a non-irreducible primary monomial ideal. We have $\height(I) = 3$. By utilizing CoCoA \cite{CoCoA}, we obtain $\mult(R/I) = 26$. Furthermore,

\[
\mult(R/I^2) = 112 \neq 26\binom{3+2-1}{2-1} = 104;
\]

\[
\mult(R/I^3) = 294 \neq 26\binom{3+3-1}{3-1} = 260;
\]

\[
\mult(R/I^4) = 608 \neq 26\binom{3+4-1}{4-1} = 520.
\]
\end{Example}
\vspace{2mm}

{\bf \noindent Acknowledgment:}
The first author acknowledges the assistance of Jiawei Bao in computations. Many of the computations related to this project were done using CoCoA \cite{CoCoA}.

\vspace{2mm}

{\bf\noindent Statement:} On behalf of all authors, the corresponding author states that there is
no conflict of interest.

\end{document}